\documentclass[11pt]{amsart}
\usepackage{amscd,amssymb,verbatim}
\theoremstyle{plain}
\newtheorem{Thm}{Theorem}[section]
\newtheorem{Cor}[Thm]{Corollary}
\newtheorem{Lem}[Thm]{Lemma}
\newtheorem{Prop}[Thm]{Proposition}
\newtheorem{Def}[Thm]{Definition}

\numberwithin{equation}{section}

\begin{document}
\title[Operators preserving asymptotic \(\ell_1\) spaces]
{Operators on \(C[0,1]\) preserving \\
 copies of asymptotic \(\ell_1\) spaces}
\author{I. Gasparis}
\address{Department of Mathematics \\
Aristotle University of Thessaloniki \\
Thessaloniki 54124, Greece.}
\email{ioagaspa@auth.gr}
\begin{abstract}
Given separable Banach spaces \(X\), \(Y\), \(Z\) and a bounded linear operator
\(T \colon X \to Y\), then \(T\) is said to preserve a copy of \(Z\) 
provided that there exists a closed linear subspace \(E\) of \(X\) isomorphic
to \(Z\) and such that the restriction of \(T\) to \(E\) is an into isomorphism.
It is proved that every operator on \(C([0,1])\) which preserves a copy of an asymptotic
\(\ell_1\) space also preserves a copy of \(C([0,1])\). 
\end{abstract}
\maketitle
\section{Introduction}
The complementation problem for the Banach space \(C(K)\), \(K\) compact metrizable, asserts that every 
complemented subspace \(E\) of \(C(K)\) is isomorphic to \(C(L)\) for some compact metric space \(L\)
(\cite{bp3}, \cite{p2}, \cite{lw}). For an in-depth analysis of this problem, we refer to \cite{r3}.
We recall that a subspace \(Y\) of a Banach space \(X\) is {\em complemented}
if \(Y\) is the range of an idempotent operator on \(X\). (In the sequel, by a subspace
of a Banach space we shall always mean an infinite-dimensional, closed, linear subspace. 
All operators will be assumed to be bounded and linear.)
Only partial results are known regarding this complementation problem. H. Rosenthal \cite{r1}
showed that if \(E^*\) is non-separable, then \(E\) is isomorphic to \(C([0,1])\). When \(E^*\)
is separable, then the results of D. Alspach and Y. Benyamini (\cite{ab}, \cite{be}) 
yield a countable compact metric space \(L\) such that each one of \(E\) and \(C(L)\) is
isomorphic to a quotient of the other. The \(L\) in question can be determined by an ordinal
index called the Szlenk index of \(E\) \cite{sz}. 

There is lack of sufficient understanding about
what a projection operator on \(C(K)\) is and, in fact, all known results about complemented subspaces
of \(C(K)\) follow from results about general operators on \(C(K)\). Following \cite{fgj}, given
separable Banach spaces \(X\), \(Y\) and \(Z\) then an operator \(T \colon X \to Y\)
is said to {\em preserve a copy} of \(Z\), if there exists a subspace \(Z'\) of \(X\)
isomorphic to \(Z\) and such that the restriction of \(T\) to \(Z'\) is an into isomorphism.
It is a classical result, due to A. Pelczynski \cite{p2}, that a non-weakly compact operator defined
on a \(C(K)\) space, preserves a copy of \(c_0\).
Rosenthal \cite{r1} actually proved that every operator \(T \colon C(K) \to X \), where \(X\) is separable,
for which \(T^*(X^*)\) is non-separable, preserves a copy of \(C(K)\). This result, combined with those
of Pelczynski (\cite{p2}, \cite{p3}) and Miljutin \cite{mi}, yields Rosenthal's result about complemented
subspaces of \(C(K)\). 
In order to extend Rosenthal's result to operators \(T\) such that
\(T^*(X^*)\) is separable, one needs to determine
the largest countable ordinal \(\alpha\) such that \(T\) preserves a copy of \(C([1, \alpha])\) and show
that \(\alpha\) is comparable to the Szlenk index of \(T^*(B_{X^*})\) (\cite{r3}, \cite{sz}).
D. Alspach \cite{a2} and J. Bourgain \cite{b2} showed that such an extension is possible when the Szlenk index 
is sufficiently large. However, such an extension fails to exist in general (\cite{a2}, \cite{a3}).

We note that it is unknown whether or not a complemented subspace \(E\) of \(C(K)\) with \(E^*\) separable is
{\em \(c_0\)-saturated}, that is every subspace of \(E\) contains an isomorph of \(c_0\). The latter
is a well known property of \(C(K)\) spaces for countable \(K\) \cite{ps}.
Regarding this problem, it is important to determine conditions on an operator
\(T \colon C(K) \to X\), \(X\) separable, which ensure that \(T^*(X^*)\) is non-separable.
Rosenthal conjectured that this is the case when \(T\) preserves a copy of
a Banach space not containing an isomorph of \(c_0\). 
In particular, if \(T\) preserves a copy of a reflexive space, then the conjecture
asserts that \(T\) preserves 
a copy of the universal space \(C([0,1])\).
An affirmative answer to Rosenthal's conjecture
yields immediately that every complemented subspace of \(C(K)\) with separable dual is \(c_0\)-saturated.
Bourgain \cite{b} showed that Rosenthal's conjecture is valid when
\(T\) preserves a copy of a Banach space \(Y\) having non-trivial cotype. 
In particular, if \(T\) preserves a copy of some \(\ell_p\), \(1 \leq p < \infty\), then
\(T\) preserves a copy of \(C(K)\) (the case \(p=1\) follows directly from \cite{r1}).
By the results of \cite{mp}, \(Y\) has non-trivial cotype if, and only if, \(c_0\) is not
finitely representable in \(Y\). This property is stronger than that of the non-embeddability of
\(c_0\) into \(Y\). Indeed, there exist reflexive spaces for which \(c_0\) is finitely representable
in every subspace. Schlumprecht's space \cite{sch} has this property, as is shown by
D. Kutzarova and P.K. Lin in \cite{kl}. A 
sub-class of the class of mixed Tsirelson spaces (defined in \cite{ad1}), was shown by
S. Argyros, I. Deliyanni, D. Kutzarova and A. Manoussakis in \cite{ad} to share this property as well.
E. Odell and Th. Schlumprecht \cite{o4} constructed a reflexive space with an unconditional basis
for which \(\ell_p\) (\(1 \leq p < \infty\)) and \(c_0\) are block-finitely represented in all
block bases. They have also constructed \cite{o5} a reflexive space with a basis \((e_n)\) such that every
finite, monotone basic sequence is block-finitely represented in every block basis of \((e_n)\).

It thus seems natural to study operators on \(C(K)\) spaces which preseve copies of
Banach spaces not containing \(c_0\), yet \(c_0\) is finitely representable in every subspace,
and attempt to verify Rosenthal's conjecture for this class of spaces.
The present paper aims towards this direction. We consider spaces from an important family
of Banach spaces that includes 
the mixed Tsirelson spaces of \cite{ad} cited above. These are the asymptotic \(\ell_1\)
spaces (\cite{mmt}, \cite{mt}). 
We recall here that a Banach space is called {\em asymptotic} \(\ell_1\) with respect to its Schauder basis
\((e_n)\) \cite{mt}, if there exists a constant \(C > 0 \) such that for every block basis
\(u_1 < \dots < u_m \) of \((e_n)\) with \(m \leq \min \mathrm{ supp } \, u_1\), we have that
\[ \| \sum_{i=1}^m u_i \| \geq C \sum_{i=1}^m \|u_i \|.\]
For an in-depth study of asymptotic \(\ell_1\) spaces, we refer to
\cite{mmt}, \cite{mt}, \cite{otw}, \cite{o2}.

Our main result is as follows:
\begin{Thm} \label{T1}
Let \(T \colon C(K) \to C(K)\) be an operator preserving a copy of an asymptotic \(\ell_1\)
space. Then \(T\) preserves a copy of \(C(K)\).
\end{Thm}
Note that \(K\) is necessarily uncountable and therefore, in view of Miljutin's theorem,
we shall assume that \(K\) is a totally disconnected, uncountable, compact metrizable space.
In particular, we may take \(K\) in the statement of Theorem \ref{T1} to be the Cantor discontinuum.
Theorem \ref{T1} is a direct consequence of Corollary \ref{C3}, proved in Section
\ref{S4}, and Rosenthal's result \cite{r1}. We show that if \(\mathcal{M}\) is a 
\(w^*\)-compact subset of \(B_{C(K)^{*}}\) which norms a subspace \(X \subset C(K)\)
which is asymptotic \(\ell_1\),
then \(\mathcal{M}\) is not separable in norm. Recall that \(\mathcal{M}\) is said to
{\em norm } \(X \), if there exists a constant \(\rho > 0\)
such that \(|\int_K \, f \, d \mu | \geq \rho \|f\|\), for all \(f \in X\).

Our proof of Theorem \ref{T1} depends on ideas and results from \cite{b}. In Section \ref{S2}
we give a criterion for a \(w^*\)-compact set \(\mathcal{M}\)
consisting of positive measures on \(K\), to be non-separable in the \(C(K)^*\)-norm.
More precisely, we show the following:
\begin{Thm} \label{T2}
Let \(\mathcal{M} \subset B_{C(K)^*}\) be a \(w^*\)-compact set of positive measures on \(K\) and let
\((f_n)\) be a normalized weakly null sequence in \(C(K)\).
Suppose that there exist a scalar \(\rho > 0\), a sequence of
positive scalars \((\epsilon_n)\) and \(N \in [\mathbb{N}]\)
such that for every \(L \in [N]\), \(L= (l_i)\), and all \(n \in \mathbb{N}\)
there exists some \(\mu \in \mathcal{M}\) satisfying
\(\mu([|f_{l_{2i}}| \geq \epsilon_{l_{2i-1}}]) \geq \rho\), for all \(i \leq n\).
Then \(\mathcal{M}\) is not separable in norm.
\end{Thm}
Note that Lemma 4 of \cite{b}, allows us to reduce the proof of Corollary \ref{C3}
to the case of \(w^*\)-compact sets of positive measures and then apply the preceding result
to conclude non-separability. This approach is different from that of \cite{b}, where non-separability
in norm follows by showing that the Szlenk index of \(\mathcal{M}\) is equal to \(\omega_1\), the first
uncountable ordinal.

A crucial step in the proof of Bourgain's result is an inequality, Lemma 6 in \cite{b}, which
holds in spaces with non-trivial cotype. In our case, the lack of cotype is substituted by a property
of asymptotic \(\ell_1\) spaces, Theorem \ref{T3}, which, loosely speaking, asserts that in an asymptotic
\(\ell_1\) subspace \(X\) of \(C(K)\) having a normalized basis of non-negative functions \((f_n)\),
one can find, for all \(\alpha < \omega_1\), a normalized block 
\(u = \sum_n a_n f_n \) and a \(t_0 \in K\) such that
\((u | I)(t_0)\) essentially estimates the norm of \( u | I\), for all \(I \subset \mathrm{ supp } \, u\)
for which \(\|u | I\|\) is significant, and moreover, \(\|u | I \|\) is negligible whenever
\(I\) belongs to the \(\alpha\)-th Schreier class \cite{aa}. 
Theorem \ref{T3} is applied in Proposition \ref{mP} which is the main step for proving Corollary \ref{C3}.
\section{Preliminaries} \label{S1}
We shall make use of standard Banach space facts and terminology
as may be found in \cite{lt}. All Banach spaces considered in this paper
are real, infinite-dimensional. For a Banach space \(E\) we let \(B_E\)
denote its closed unit ball.

If \(X\) is any set, we let
\([X]^{< \infty}\) denote the set of its finite subsets,
while \([X]\) stands for the set of all infinite subsets
of \(X\). If \(M \in [\mathbb{N}]\), we shall adopt the convenient
notation \(M=(m_i)\) to denote the increasing enumeration of
the elements of \(M\).

A family \(\mathcal{F} \subset [\mathbb{N}]^{< \infty}\) is {\em
hereditary} if \(G \in \mathcal{F}\) whenever \(G \subset F\) and
\(F \in \mathcal{F}\). \(\mathcal{F}\) is {\em compact}, if it is compact
with respect to the topology of pointwise convergence in
\([\mathbb{N}]^{< \infty}\).

If \(E\) and \(F\) are finite subsets of
\(\mathbb{N}\), we write \(E < F\) when
\(\max E < \min F\).

{\bf Notation}. Given \(\mathcal{F} \subset [\mathbb{N}]^{< \infty}\)
and \(M \in [\mathbb{N}]\), we set
\(\mathcal{F}[M] = \{ F \cap M : \, F \in \mathcal{F} \}\).
Clearly, \(\mathcal{F}[M]\) is hereditary (resp. compact), if
\(\mathcal{F}\) is.

We shall now recall the transfinite definition of
the Schreier families \(S_\xi\), \(\xi < \omega_1\) \cite{aa}.
First, given a countable ordinal \(\alpha\) we associate
to it a sequence of successor ordinals, \((\alpha_n +1)\),
in the following manner:
If \(\alpha\) is a successor ordinal we let
\(\alpha_n = \alpha - 1\) for all \(n\).
In case \(\alpha\) is a limit ordinal, we choose
\((\alpha_n +1)\) to be a strictly increasing sequence of ordinals
tending to \(\alpha\).

Now set \(S_0 = \{\{n\}: \, n \in \mathbb{N}\} \cup \{\emptyset\}\)
and
\(S_1 = \{F \subset \mathbb{N}: \, |F| \leq \min F\}
\cup \{\emptyset\}\). 
Let \(\xi < \omega_1\) and assume \(S_\alpha\) has been defined
for all \(\alpha < \xi\). If \(\xi\) is a successor
ordinal, say
\(\xi = \zeta + 1\), define
\[S_{\xi } = \{ \cup_{i=1}^n F_i : \, n \in \mathbb{N}, \, F_i \in S_\zeta \, 
\forall \, i \leq n, \,  
 n \leq \min F_1, \, F_1 < \dots < F_n \}.\]
In the case \(\xi\) is a limit ordinal,
let \((\xi_n + 1)\) be the
sequence of successor ordinals associated to \(\xi\).
Set
\[S_\xi = \cup_n \{ F \in S_{\xi_n + 1}: \, n \leq \min F\}
\cup \{\emptyset\}.\]
It is shown in \cite{aa} that the Schreier
family \(S_\xi\) is hereditary and compact
for all \(\xi < \omega_1\). 
The Schreier families played an important role in the recent development of Banach
space theory. For a detailed exposition of this development and the use of ordinal indices
in Banach space theory, we refer to \cite{o2}, \cite{o3}.
 
Given a sequence \((e_n)\) of non-zero vectors in some Banach space \(X\), then
the vector \(u \in X\) is called a {\em block} of \((e_n)\) if 
\(u = \sum_{i \in I} \lambda_i e_i\) where \(I \in [\mathbb{N}]^{< \infty}\)
and \((\lambda_i)_{i \in I}\) are scalars. We also denote by \(u | J\) the vector
\(\sum_{i \in I \cap J} \lambda_i e_i\) for all \(J \subset \mathbb{N}\). 

A (finite or infinite) sequence of non-zero vectors \((u_n)\) in \(X\) is called
a {\em block subsequence} of \((e_n)\), if there exists a sequence of non-zero scalars \((\lambda_n)\),
and a sequence \((F_n)\) of finite subsets of \(\mathbb{N}\) with \(F_1 < F_2 < \dots \),
such that \(u_n = \sum_{i \in F_n} \lambda_i e_i\), for all \(n \in \mathbb{N}\).
We then call \(F_n\), the {\em support} of \(u_n\) and write
\(F_n = \mathrm{ supp } \, u_n \) for all \(n \in \mathbb{N}\). The notation 
\(u_1 < u_2 < \dots \) indicates that \(F_1 < F_2 < \dots \). 
In case \((e_n)\) is a {\em basic} sequence, that is \((e_n)\) is a Schauder basis for
its closed linear span, then we call \((u_n)\) a {\em block basis} of \((e_n)\).
\section{A criterion for the norm-separability of \(w^*\)-compact sets of positive measures} \label{S2}
In this section we give the proof of Theorem \ref{T2}.
This result will follow after establishing the next
\begin{Prop} \label{P1}
Let \(\mathcal{M} \subset B_{C(K)^*}\) be a \(w^*\)-compact set of positive measures on \(K\) and let
\((f_n)\) be a normalized weakly null sequence in \(C(K)\).
Let \(\rho \) be a positive scalar and let \((\epsilon_n)\), \((\delta_n)\) be sequences of
positive scalars. 
Then there exists a subsequence \((f_{m_n})\) of \((f_n)\) satisfying the following property:
whenever \(I \in [\mathbb{N}]^{< \infty}\) and \(\mu \in \mathcal{M}\)
are so that \(\mu([|f_{m_{2i}}| \geq \epsilon_{m_{2i-1}}]) \geq \rho\),
for all \(i \in I\), there exists some \(\nu \in \mathcal{M}\) such that
\(\nu([|f_{m_{2i}}| \geq \epsilon_{m_{2i-1}}]) \geq \rho\) for all \(i \in I\),
yet, \(\nu(|f_{m_{2i}}|) \leq \delta_{m_{2i-1}}\) for all \(i \in \mathbb{N} \setminus I\). 
\end{Prop}
The technique for proving this proposition is based on the infinite Ramsey theorem (\cite{ell}, \cite{o1}),
and is similar to methods developed in \cite{e}, \cite{o1}, \cite{amt}, \cite{ag}
for the study of subsequential properties of weakly null sequences.
We postpone the proof of Proposition \ref{P1} in order to give the
\begin{proof}[Proof of Theorem \ref{T2}]
The hypotheses of the theorem readily imply that \(\lim_{n \in N} \epsilon_n = 0\).
Thus, there will be no loss of generality in assuming that \(\epsilon_n < \rho /2\),
for all \(n \in N\).
Write \(N = (k_n)\) and apply Proposition \ref{P1} to the weakly null sequence 
\((f_{k_n})\), the set of measures \(\mathcal{M}\), the scalar
\(\rho\) and the scalar sequences "\((\epsilon_n)\)"\(=(\epsilon_{k_n})\) and 
"\((\delta_n)\)"\(=(\epsilon_{k_n}^2)\). Taking in account the hypothesis of the theorem,
we shall obtain \(P \in [N]\), \(P = (p_n)\), with the following property: for every
\(I \in [\mathbb{N}]^{< \infty}\) there exists some \(\mu \in \mathcal{M}\) such that
\[\mu([|f_{p_{2i}}| \geq \epsilon_{p_{2i-1}}]) \geq \rho \text{ for all } i \in I,
\text{ yet, } \mu(|f_{p_{2i}}|) \leq \epsilon_{p_{2i-1}}^2 \text{ for all } i \in \mathbb{N} \setminus I.\] 
Now let \(L \in [\mathbb{N}]\), \(L=(l_n)\). Our previous work yields, for every 
\(n \in \mathbb{N}\), some \(\mu_n \in \mathcal{M}\) satisfying 
\begin{align}
&\mu_n([|f_{p_{2l_i}}| \geq \epsilon_{p_{2l_i-1}}]) \geq \rho \text{ for all } i \leq n, 
\text{ yet, } \notag \\
&\mu_n(|f_{p_{2j}}|) \leq \epsilon_{p_{2j-1}}^2 \text{ for all } j \in \mathbb{N} \setminus 
\{l_1, \dots , l_n\}. \notag
\end{align} 
Letting \(\mu_L \in \mathcal{M}\) be any \(w^*\)-cluster point of \((\mu_n)\), we infer from the above that
\[\mu_L([|f_{p_{2j}}| \geq \epsilon_{p_{2j-1}}]) \geq \rho \text{ for all } j \in L,
\, \mu_L(|f_{p_{2j}}|) \leq \epsilon_{p_{2j-1}}^2 \text{ for all } j \in \mathbb{N} \setminus L.\] 
Note that in obtaining the first inequality above, we used the fact that \(\tau(F) \geq \delta\)
whenever \(F\) is a closed subset of \(K\) and \(\tau\) is the \(w^*\)-limit in \(C(K)^*\) of 
a sequence of positive measures \((\tau_n)\) satisfying \(\tau_n(F) \geq \delta > 0\) for
all \(n \in \mathbb{N}\).

We obtain in particular, that 
\(\mu_L([|f_{p_{2j}}| \geq \epsilon_{p_{2j-1}}]) \leq \epsilon_{p_{2j-1}}\), for all
\(j \in \mathbb{N} \setminus L\).

Finally, if \(L_1\) and \(L_2\) are distinct members of \([\mathbb{N}]\) we may choose
without loss of generality some \(j \in L_1 \setminus L_2\). We now have, by the manner the
\(\mu_L\)'s have been selected, that
\[\mu_{L_1}([|f_{p_{2j}}| \geq \epsilon_{p_{2j-1}}]) \geq \rho, \text{ while }  
\mu_{L_2}([|f_{p_{2j}}| \geq \epsilon_{p_{2j-1}}]) \leq \epsilon_{p_{2j-1}}.\]
We deduce from the above, that \(\|\mu_{L_1} - \mu_{L_2} \| \geq \rho - \epsilon_{p_{2j-1}} \geq
\rho / 2 \).
Therefore, \(\mathcal{M}\) is not separable in norm, as claimed.
\end{proof}
Before giving the proof of Proposition \ref{P1}, we need to introduce some notation and
terminology. We first fix a normalized weakly null sequence \((f_n)\) in \(C(K)\), a \(w^*\)-compact
subset \(\mathcal{M}\) of \(B_{C(K)^*}\) consisting of positive measures on \(K\), a positive
scalar \(\rho\) and sequences of positive scalars \((\epsilon_n)\) and \((\delta_n)\).

{\bf Notation 1}. For a finite subset \(A\) of \(\mathbb{N}\) of even cardinality (including the case
\(A = \emptyset\)), say \(A= \{m_1, < \dots , < m_{2n}\}\), we set
\(A^{(2)} = \{m_{2i} : i \leq n \}\). Given \(m \in A^{(2)}\), \(m = m_{2i}\) 
for some \(i \leq n\), we set \(m^{-} = m_{2i-1}\). If \(F \subset A^{(2)}\)
we set \(F^{-} = \{ m^{-} : m \in F \}\).

{\bf Notation 2}. Let \(A \in [\mathbb{N}]^{< \infty}\) be of even cardinality and
let \(L \in [\mathbb{N}]\), \(L = (l_i)\), with \(\max A < l_1 \) (\(\max \emptyset = 0\)).
Given \(F \subset A^{(2)}\) we set
\[ (F,A) \sqcup L = F \cup F^{-} \cup \{l_j : \, j \geq 3\}.\]

{\bf Terminology 1}. A pair \((F_1, F_2)\) of finite subsets of \(\mathbb{N}\) is said
to be {\em appropriate}, if \(F_2\) is of even cardinality and \(F_1 \subset F_2^{(2)}\).

{\bf Terminology 2}. Given \(L \in [\mathbb{N}]\), \(L = (l_i)\), and \(n \in \mathbb{N}\),
then a measure \(\mu \in \mathcal{M}\) is called \((L,n)\)-{\em good} if
\(\mu([|f_{l_{2i}}| \geq \epsilon_{l_{2i-1}}]) \geq \rho\), for all \(i \leq n\).

{\bf Terminology 3}. Let \((F_1,F_2)\) be an appropriate pair and let
\(L \in [\mathbb{N}]\), \(L=(l_i)\), with \(\max F_2 < l_1 \).
\(L\) is said to be \((F_1,F_2)\)-{\em admissible}, provided the following condition
is fulfilled: If \(\mu \in \mathcal{M}\) and \(n \in \mathbb{N}\) satisfy
\[\mu \text{ is } \bigl ( (F_1, F_2) \sqcup L, n \bigr )- \text{good},\]
then there exists some
\(\nu \in \mathcal{M}\) which is 
\begin{align}
&\bigl ( (F_1, F_2) \sqcup L, n \bigr )- \text{good and such that} \notag \\ 
&\nu(|f_m|) < \delta_{m^{-}} \text{ for all } m \in F_2^{(2)} \setminus F_1,
\text{ and } \nu ( |f_{l_2}|) < \delta_{l_1}. \notag
\end{align}
Finally, \(L\) is called \(F_2\)-{\em admissible} if it is
\((F, F_2)\)-admissible for every \(F \subset F_2^{(2)}\).
\begin{proof}[Proof of Proposition \ref{P1}]
We shall inductively construct a decreasing sequence
\(M_0 \supset M_1 \supset \dots\) of infinite subsets of \(\mathbb{N}\)
and an increasing sequence \(m_1 < m_2 < \dots \) of positive integers
such that
\begin{enumerate}
\item \(m_{2i-1}\) and \(m_{2i}\) are the third and fourth elements, respectively,
of \(M_{i-1}\) for all \(i \in \mathbb{N}\).
\item Every \(L \in [M_i]\) is \(\{m_1, \dots , m_{2i}\}\)-admissible 
for all \(i \in \mathbb{N} \cup \{0\}\)
(when \(i=0\), \(\{m_1, \dots, m_{2i} \} = \emptyset\)).
\end{enumerate}
The construction of \(M_0\) is implicit in the general inductive step and
therefore it will not be discussed. So we shall assume that \(i \geq 1\)
and that \(M_0 \supset \dots \supset M_{i-1}\)  and \(m_1 < \dots < m_{2i-2}\)
have been constructed satisfying \((1)\) and \((2)\) for all \(j \leq i-1\).
Let \(m_{2i-1}\) and \(m_{2i}\) be the third and fourth element, respectively,
of \(M_{i-1}\). To simplify our notation, we set \(F_j = \{m_1, \dots, m_{2j}\}\),
for all \(j \in \mathbb{N} \cup \{0\}\). Given
\(F \subset F_i^{(2)}\) and \(Q \in [M_{i-1}]\), we define
\[\mathcal{D}_F = \{ L \in [Q]: \, m_{2i} < \min L \text{ and }
L \text{ is } (F, F_i)-\text{admissible} \}.\]
It follows directly from the definitions (see Terminologies \(2\) and \(3\)),
that \(\mathcal{D}_F\) is closed in the topology of pointwise convergence.
The infinite Ramsey theorem yields some \(P \in [Q]\) such that either
\(\mathcal{D}_F \subset [P]\), or \(\mathcal{D}_F \cap [P] = \emptyset \).
We are going to show that only the first alternative can hold. Assuming this
is not the case we obtain a contradiction as follows: Write \(P = (p_j)\)
and let \(k \in \mathbb{N}\), \(k \geq 2\). We define sets
\[L_j = \{p_1, p_j\} \cup \{p_t : \, t > k \}, \text{ for all } 2 \leq j \leq k.\]
We also define \[R=
\begin{cases}
\{q_{2i-1}, q_{2i}, m_{2i-1}, m_{2i} \} \cup \{p_t: \, t > k \}, 
&\text{ if } m_{2i} \in F ; \\
\{ m_{2i-1}, m_{2i} \} \cup \{p_t: \, t > k \}, 
&\text{ if } m_{2i} \notin F.
\end{cases}\]
In the above, we have let \(q_{2i-1}\) and \(q_{2i}\) denote the first and second
elements, respectively, of \(M_{i-1}\).
One checks that
\begin{equation} \label{e1}
(F,F_i) \sqcup L_j = (F\setminus \{m_{2i}\}, F_{i-1}) \sqcup R,
\text{ for all } 2 \leq j \leq k.
\end{equation}
Since \(L_j \in [P]\) we have that \(L_j\) is not \((F, F_i)\)-admissible
for all \(2 \leq j \leq k\). It follows that for every \(2 \leq j \leq k\)
there exist \(\mu_j \in \mathcal{M}\) and \(d_j \in \mathbb{N}\) so that
\begin{equation} \label{e2}
\mu_j \text{ is } \bigl ( (F,F_i) \sqcup L_j, d_j \bigr )-\text{good}
\end{equation}
and moreover,
\begin{align} \label{e3}
&\text{every } \nu \in \mathcal{M} \text{ which is }
\bigl ( (F,F_i) \sqcup L_j, d_j \bigr )-\text{good}  \\
&\text{ and satisfies } \nu(|f_m|) < \delta_{m^{-}}, \text{ for all }
m \in F_i^{(2)} \setminus F, \notag \\
&\text{ also satisfies } \nu(|f_{p_j}|) \geq \delta_{p_1}. \notag
\end{align}
Next, choose \(j_0 \leq k \) such that \(d_{j_0} = \max \{d_j : \, 2 \leq j \leq k \}\).
We infer from \eqref{e2} and \eqref{e1}, that \(\mu_{j_0}\) is 
\(\bigl ( (F\setminus \{m_{2i}\}, F_{i-1}) \sqcup R , d_{j_0} \bigr ) \)-good. Since \(R \in [M_{i-1}]\),
it is \((F\setminus \{m_{2i}\}, F_{i-1})\)-admissible, by the inductive hypothesis.
Hence there exists some \(\nu \in \mathcal{M}\) satisfying
\begin{align} \label{e4}
&\nu \text{ is }
\bigl ( (F \setminus \{m_{2i}\} ,F_{i-1}) \sqcup R, d_{j_0} \bigr )-\text{good}, \\  
&\nu(|f_m|) < \delta_{m^{-}}, \text{ for all }
m \in F_{i-1}^{(2)} \setminus (F \setminus \{m_{2i}\}), \notag \\
&\nu(|f_{r_2}|) < \delta_{r_1}. \notag
\end{align}
In the above, \(r_1 < r_2 \) are the first two elements of \(R\).
We now observe that the definition of \(R\) and \eqref{e4} lead to 
\begin{equation} \label{e5}
\nu( (|f_m|) < \delta_{m^{-}}, \text{ for all }
m \in F_i^{(2)} \setminus F.
\end{equation}
Note also that since \(d_j \leq d_{j_0}\) for all \( 2 \leq j \leq k\),
\eqref{e4} and \eqref{e1} yield that 
\begin{equation} \label{e6}
\nu \text{ is } \bigl ( (F,F_i) \sqcup L_j, d_j \bigr )-\text{good, for all }  
2 \leq j \leq k.
\end{equation}
By combining \eqref{e6} with \eqref{e5} and \eqref{e3}, we conclude that
\[\nu (|f_{p_j}|) \geq \delta_{p_1}, \text{ for all }
2 \leq j \leq k.\]
Since \(k \geq 2\) was arbitrary, this contradicts the assumption that
\((f_n)\) is weakly null.

Therefore, we have indeed that \(\mathcal{D}_F \subset [P] \) and thus
every infinite subset of \(P\) is \((F,F_i)\)-admissible.
If we now let \(\{H_1, \dots ,H_s\}\) be an enumeration 
of the subsets of \(F_i^{(2)}\), then successive applications of the preceding
argument yield infinite subsets \(P_1 \supset \dots \supset P_s \)
of \(M_{i-1}\) such that \(L\) is \((H_j , F_i)\)-admissible for
all \(L \in [P_j]\) and every \(j \leq s\). It follows now that 
\(M_i=P_s\) has the property that every \(L \in [M_i]\) is
\(F_i\)-admissible. This completes the inductive construction.

We finally show that the subsequence \((f_{m_n})\) of \((f_n)\) satisfies the
conclusion of the proposition. Indeed, let 
\(I \in [\mathbb{N}]^{< \infty}\) and \(\mu \in \mathcal{M}\)
satisfy \(\mu([|f_{m_{2i}}| \geq \epsilon_{m_{2i-1}}]) \geq \rho\),
for all \(i \in I\). Set \(p = \max I\). It follows that for every
integer \(n \geq p\), \(\mu\) is 
\[\bigl ( (\{m_{2i} : i \in I \}, \{m_1, \dots , m_{2n} \} )
\sqcup M_n, |I| \bigr )-\text{good}.\]
Since \(M_n\) is \(\{m_1, \dots , m_{2n}\}\)-admissible, for all \(n \in \mathbb{N}\),
we infer that for every \(n \geq p\), there exists some \(\mu_n \in \mathcal{M}\)
satisfying  
\begin{align}
&\mu_n \text{ is }
\bigl ( (\{m_{2i} : i \in I \}, \{m_1, \dots , m_{2n} \} )
\sqcup M_n, |I| \bigr )-\text{good}, \notag \\
&\mu_n(|f_{m_{2i}}|) < \delta_{m_{2i-1}}, \text{ for all }
i \in \{1, \dots, n \} \setminus I. \notag
\end{align}
Hence, for all \(n \geq p\) there exists \(\mu_n \in \mathcal{M}\)
satisfying
\begin{align}
&\mu_n([|f_{m_{2i}}| \geq \epsilon_{m_{2i-1}}]) \geq \rho,
\text{ for all } i \in I, \notag \\
&\mu_n(|f_{m_{2i}}|) < \delta_{m_{2i-1}}, \text{ for all }
i \in \{1, \dots, n \} \setminus I. \notag
\end{align}
Finally, let \(\nu \in \mathcal{M}\) be a \(w^*\)-cluster point of
the sequence \((\mu_n)_{n \geq p}\). Clearly, \(\nu\) is as desired.
\end{proof}
\section{A property of asymptotic \(\ell_1\) spaces} \label{S3}
We start this section with a simple observation about the usual \(\ell_1\)-basis \((e_n)\):
If \(u = (1/n) \sum_{i=1}^n e_i\) and \(u^* = \sum_{i=1}^n e_i^* \in B_{\ell_\infty}\), then
\(u^*( u | F) = \|u | F \|\) for every \(F \subset \mathrm{ supp } \, u\). Moreover,
\(\|u | \{i\}\| = 1/n < \epsilon\), if \(n\) is sufficiently large, for all \(i \leq n\).
In this section, we investigate if a similar property holds in asymptotic
\(\ell_1\) spaces. This is the content of Theorem \ref{T3} below.
In what follows, \(K\) is a compact metrizable space and \(X\) is a subspace of \(C(K)\)
which is \(C\)-asymptotic \(\ell_1\) with respect to its normalized Schauder basis \((e_n)\).
This means that there exists a constant \(C> 0\) such that
every block basis \(u_1 < \dots < u_m\) of \((e_n)\) with 
\(m \leq \min \mathrm{ supp } \, u_1\) satisfies
\(\|\sum_{i=1}^m u_i \| \geq C \sum_{i=1}^m \|u_i\|\).

We shall need to work with non-negative functions in \(C(K)\) and so we note that if
\((x_n)\) is a basic sequence in \(C(K)\)
then its sequence of absolute values, \((|x_n|)\), may not belong to the closed linear span
of \((x_n)\). However, as is shown in Lemma 5 of \cite{b}, under certain conditions, \((|x_n|)\)
may somehow inherit properties of \((x_n)\). We make this more precise below, by adapting the aforementioned
lemma of \cite{b} into the context of asymptotic \(\ell_1\) spaces.
\begin{Lem} \label{Lex}
We set \(f_n = |e_n|\) for all \(n \in \mathbb{N}\). Then 
for every block subsequence \(u_1 < \dots < u_m\) of \((f_n)\) with \(m \leq \min \mathrm{ supp } \, u_1\)
and such that each \(u_i\) is a positive linear combination of the \(f_n\)'s, one has
that \(\|\sum_{i=1}^m u_i \| \geq C \sum_{i=1}^m \|u_i\|\).
\end{Lem}
\begin{proof}
Clearly, \((f_n)\) is
a normalized sequence of non-negative functions in \(C(K)\). 
Write \(u_i = \sum_{j \in F_i} a_j f_j\) where \(a_j > 0\) for all
\(j \in F_i\) and all \(i \leq m\). Of course, \(F_1 < \dots < F_m\). We next choose, for every \(i \leq m\),
\(t_i \in K\) and a sign \(\sigma_j \), for all \(j \in F_i\), so that
\(\|u_i \| = \sum_{j \in F_i} a_j \sigma_j e_j(t_i)\). Put \(v_i = \sum_{j \in F_i} a_j \sigma_j e_j\)
for all \(i \leq m\). It is easily checked that \(\|v_i \| = \|u_i \|\) for all \(i \leq m\), and that
\(\|\sum_{i=1}^m v_i \| \leq \|\sum_{i=1}^m u_i\|\). Since \(\mathrm{ supp } \, v_i = F_i\), 
\(v_1 < \dots < v_m \) is a block basis of \((e_n)\) with \(m \leq \mathrm{ supp } \, v_1\).
The above, clearly prove the assertion as \(X\) is \(C\)-asymptotic \(\ell_1\) with respect to
\((e_n)\).
\end{proof}
\begin{Def}
Given \(0 < \epsilon < 1\) and \(\alpha < \omega_1\), then a normalized block
\(u = \sum_{i \in I} \lambda_i f_i\) of \((f_n)\) with \(I \in [\mathbb{N}]^{< \infty}\) 
and \(\lambda_i > 0\) for all \(i \in I\), is called an \((\alpha, \epsilon)\) block
provided the following conditions hold:
\begin{enumerate}
\item There exists a \(t \in K\) such that for every \(J \subset I\) with
\(\|\sum_{i \in J} \lambda_i f_i \| \geq \epsilon \) we have that
\(\|\sum_{i \in J} \lambda_i f_i \| \leq (1 + \epsilon) \sum_{i \in J} \lambda_i f_i(t)\).
\item \(\|\sum_{i \in J} \lambda_i f_i \| < \epsilon^2\), for all \(J \subset I\) with
\(J \in S_\alpha\).
\end{enumerate}
Any \(t \in K\) satisfying \((1)\) will be said to strongly norm the \((\alpha, \epsilon)\) block \(u\). 
\end{Def}
The main result of this section is the following
\begin{Thm} \label{T3}
Let \(K\) be a compact metrizable space and let \(X\) be a subspace
of \(C(K)\) which is asymptotic \(\ell_1\) with respect to its normalized
Schauder basis \((e_n)\). Set \(f_n = |e_n|\) for all \(n \in \mathbb{N}\). Then
for every \(0 < \epsilon < 1\) and \(\alpha < \omega_1\) and all \(N \in [\mathbb{N}]\) there
exists an \((\alpha, \epsilon)\) block of \((f_n)\) supported by \(N\).
\end{Thm}
The proof of Theorem \ref{T3}, requires a few intermediate steps which are presented below.
The first one is a simple permanence property of \((\alpha, \epsilon)\) blocks.
\begin{Lem} \label{L1}
Suppose that \(u\) is an \((\alpha, \epsilon)\) block strongly normed by \(t_0 \in K\).
If \(I_0 \subset \mathrm{ supp } \, u\) satisfies \(\|u | I_0\| \geq \epsilon^{1/2}\),
then \(v = \frac{u | I_0}{\|u | I_0 \|}\) is an \((\alpha, \epsilon^{1/2})\) block, strongly
normed by \(t_0\) as well.
\end{Lem}
\begin{proof}
Let \(I = \mathrm{ supp } \, u\) and write \(u = \sum_{i \in I} \lambda_i f_i\),
where \(\lambda_i > 0\) for all \(i \in I\).
Let \(I_0 \subset I\) satisfy \(\|\sum_{i \in I_0} \lambda_i f_i \| \geq \epsilon^{1/2}\)
and set \(D = \|\sum_{i \in I_0} \lambda_i f_i \|\). Suppose that \(J \subset I_0\)
satisfies \(\|(1/D) \sum_{i \in J} \lambda_i f_i \| \geq \epsilon^{1/2}\).
Then, clearly, \( \| \sum_{i \in J} \lambda_i f_i \| \geq \epsilon\) and thus
\[ \| \sum_{i \in J} \lambda_i f_i \| \leq ( 1 + \epsilon) \sum_{i \in J} \lambda_i f_i(t_0),\]
as \(t_0\) strongly norms \(u\). This of course implies that
\[ \| v | J \| = \| (1/D ) \sum_{i \in J} \lambda_i f_i \| \leq ( 1 + \epsilon) 
(1/D) \sum_{i \in J} \lambda_i f_i (t_0) \leq  (1 + \epsilon^{1/2}) ( v | J)(t_0)\]
and so \(t_0\) strongly norms \(v\). Finally, suppose \(J \subset I_0\) belongs to \(S_\alpha\).
Then, \(\| \sum_{i \in J} \lambda_i f_i \| < \epsilon^2\), by our assumptions. 
We obtain that \(\|(1/D) \sum_{i \in J} \lambda_i f_i \| < \epsilon^{3/2} < \epsilon\)
which completes the proof of the lemma.
\end{proof}
Theorem \ref{T3} will be proved by transfinite induction. 
Our next lemma gives the first inductive step. 
\begin{Lem} \label{L0}
For every \(0 < \epsilon < 1\) there exists a \((0,\epsilon)\) block of
\((f_n)\) supported by \(N\).
\end{Lem}
\begin{proof}
Lemma \ref{Lex} yields the existence of some \(C > 0\) 
satisfying \(\|\sum_{i=1}^m u_i \| \geq C \sum_{i=1}^m \|u_i\|\), for every block subsequence
\(u_1 < \dots < u_m\) of \((f_n)\) with \(m \leq \min \mathrm{ supp } \, u_1\)
and such that each \(u_i\) is a positive linear combination of the \(f_n\)'s.
We define
\[ \tau = \sup \{ c > 0 : \, \forall \, M \in [N] \, \,  \exists \, L \in [M], \,
L = (l_i), \text{ with } 
\| \sum_{i=1}^{l_1} f_{l_i} \| \geq c \, l_1 \}.\]
Clearly, \(C \leq \tau \leq 1\). We remark that the modulus \(\tau\) is a special case of the
\(\delta_\alpha\)-moduli, introduced in \cite{otw} for the study of asymptotic \(\ell_1\) spaces.  
Given \(\delta > 0\), \(\delta < \tau / 2\), we can choose some \(M \in [N]\) so that
\(\| \sum_{i=1}^{l_1} f_{l_i} \| <  (\tau + \delta) l_1 \)
for all \(L \in [M]\), \(L = (l_i)\).

{\em Claim} \(1\). Let \(m \in M\) satisfy \((m-1)(\tau + 2 \delta) > m(\tau + \delta)\).
Then, for all \(t \in K\), the cardinality of the set
\(\{n \in M : \, n > m \text{ and } f_n(t) \geq \tau + 2 \delta\}\)
is smaller than \(m - 1\).

Indeed, were this claim false, we could choose \(m < n_1 < \dots < n_{m-1} \)
in \(M\) and \(t \in K\) satisfying
\(f_{n_i}(t) \geq \tau + 2 \delta \) for all \(i \leq m -1 \). It follows that
\[f_m(t) + \sum_{i=1}^{m-1} f_{n_i}(t) \geq (m-1)(\tau + 2 \delta) > m(\tau + \delta).\]
We can now choose \(L \in [M]\), \(L= (l_i)\), with \(l_1 = m\) and \(l_i = n_i\)
for \(i=2, \dots , m\). Hence, \(\| \sum_{i=1}^{l_1} f_{l_i} \| > (\tau + \delta) l_1 \),
contradicting the choice of \(M\).

{\em Claim} \(2\). Given \(m \in M\), there exist \(t_1 \in K\) and \(n_1 < \dots < n_m \) 
in \(M\) with \(m < n_1\) so that
\(f_{n_i}(t_1) \geq \tau - 2 \delta\) for all \(i \leq m\).

To prove this claim, first choose \(M_1 \in [M]\) with \(m < \delta \min M_1\).
The definition of \(\tau\) now yields some \(L \in [M_1]\), \(L = (l_i)\),
with \(\|\sum_{i=1}^{l_1} f_{l_i} \| \geq (\tau - \delta) l_1\).
We can now choose \(t_1 \in K\)
satisfying \(\sum_{i=1}^{l_1} f_{l_i} (t_1) \geq (\tau - \delta) l_1\).
The claim will follow once we establish that the set \(\{i \in \mathbb{N}: \, i \leq l_1 \text{ and }
f_{l_i}(t_1) \geq \tau - 2 \delta \}\) contains at least \(m\) elements.
Indeed, if that were not the case, we would have that
\[\sum_{i=1}^{l_1} f_{l_i}(t_1) \leq m + l_1 (\tau - 2 \delta) < l_1 (\tau - \delta),\]
contradicting the choice of \(t_1\). 

We are now ready for the proof of this lemma. Let \( 0 < \epsilon < 1 \) and choose
\(0 < \delta < \tau / 2\) and \(0 < \epsilon_0 < \epsilon\) satisfying
\[ \frac{\tau + 2 \delta}{(1 - \epsilon_0)(\tau - 2 \delta)} < 1 + \epsilon.\]
We can select \(M \in [N]\) satisfying \(\| \sum_{i=1}^{l_1} f_{l_i} \| <  (\tau + \delta) l_1 \)
for all \(L \in [M]\), \(L = (l_i)\).
We also choose \(m \in M\) with \((m-1)(\tau + 2 \delta) > m(\tau + \delta)\), 
and \(m_0 \in M\) with \(m < C \epsilon_0^2 m_0\). We apply Claim 2 in order to find \(t_0 \in K\),
\(n_1 < \dots < n_{m_0}\) in \(M\) with \(m_0 < n_1\) and such that
\begin{equation} \label{eL0}
f_{n_i}(t_0) \geq \tau - 2 \delta, \text{ for all } i \leq m_0.
\end{equation}
We are going to show that \(u = \frac{\sum_{i=1}^{m_0} f_{n_i}}{\|\sum_{i=1}^{m_0} f_{n_i}\|}\)
is a \((0, \epsilon)\) block supported by \(N\). To this end, set \(D = \|\sum_{i=1}^{m_0} f_{n_i}\|\)
and note that \(D \geq C m_0\). Let \(I \subset \{n_1, \dots , n_{m_0}\}\) satisfy
\(\|u | I \| \geq \epsilon\). Let \(t \in K\) and set
\[J_t = \{ i \in \mathbb{N}: \, i \leq m_0 , \text{ and } f_{n_i}(t) \geq \tau + 2 \delta \}.\]
Also let \(J = \{ i \in \mathbb{N}: \, i \leq m_0 \text{ and } n_i \in I \}\). We have the estimates:
\begin{align}
(u|I)(t) &=  (1/D) \sum_{i \in J \cap J_t  } f_{n_i}(t) + (1/D) \sum_{i \in J \setminus J_t} f_{n_i}(t) \notag \\
&\leq (m/D) + (1/D)(\tau + 2 \delta) |J|, \text{ since } |J_t| < m \text{ by Claim } 1, \notag \\
&\leq (m /Cm_0) + (1/D) \frac{\tau + 2 \delta}{\tau - 2 \delta} \sum_{i \in J} f_{n_i}(t_0), \text{ by }
\eqref{eL0} \notag \\
&\text{ and since } D \geq Cm_0, \notag \\
&< \epsilon_0^2 + \frac{\tau + 2 \delta}{\tau - 2 \delta} (u|I)(t_0), \text{ by the choice of } m_0, \notag \\
&\leq \epsilon_0 \|u|I\| + \frac{\tau + 2 \delta}{\tau - 2 \delta} (u|I)(t_0), \text{ since }
\|u|I\| \geq \epsilon. \notag
\end{align}
Since \(t \in K\) was arbitrary, we conclude that
\(\|u|I\| \leq (1-\epsilon_0)^{-1} \frac{\tau + 2 \delta}{\tau - 2 \delta} (u|I)(t_0)\).
Therefore, \(\|u|I\| \leq (1 + \epsilon) (u|I)(t_0)\).
By the choices made (recall that \(D \geq C m_0 > \epsilon_0^{-2}\)), we have also ensured that \(1/D < \epsilon^2\).
Thus, \(u\) is a \((0,\epsilon)\) block supported by \(N\) and strongly normed by \(t_0\). 
This completes the proof of the lemma.
\end{proof}
The general inductive step in the proof of Theorem \ref{T3}, requires the following definition.
\begin{Def}
Let \(1 \leq \alpha < \omega_1\) and denote by \((\alpha_n +1)\) the sequence of successor ordinals associated
to \(\alpha\). We also fix a decreasing sequence of positive scalars \((\epsilon_n)\) satisfying
\(\sum_n \epsilon_n^{1/2} < 1\). 
A finite normalized block subsequence \(u_1 < \dots < u_m\) of \((f_n)\) (\(m \geq 2\)), 
where each \(u_i\) is a positive linear
combination of the \(f_n\)'s, is called an \(\alpha\)-chain, provided the following properties are satisfied:
\begin{enumerate}
\item \(m = \min \mathrm{ supp } \, u_1\).
\item For every \(i \in \{2, \dots , m \}\), \(u_i\) is an \((\alpha_{d_{i-1}}, \epsilon_{d_{i-1}})\) block,
where \(d_i = \max \mathrm{ supp } \, u_i \) for all \(i \leq m\).
\item For every \(i \in \{2, \dots , m \}\), every \(j \leq d_{i-1}\) and every 
\(F \subset \mathrm{ supp } \, u_i\) with \(F \in S_{\alpha_j}\), we have that
\(F \in S_{\alpha_{d_{i-1}}}\).
\end{enumerate}
\end{Def}
We remark that the definition of an \(\alpha\)-chain is similar to that of the \((\alpha, \beta, \epsilon )\)
{\em averages} and {\em rapidly increasing sequences of special convex combinations}
(\cite{sch},\cite{gm}, \cite{ad1}, \cite{otw}). The next lemma guarantees the existence of \(\alpha\)-chains.
\begin{Lem} \label{L2}
Let \(1 \leq \alpha < \omega_1\) and assume that the conclusion of Theorem \ref{T3}
holds for all \(\beta < \alpha\). Given \(N \in [\mathbb{N}]\) there exists \(N_0 \in [N]\)
such that setting 
\begin{align}
\tau = \sup \{ &c > 0 : \, \forall \, L \in [N_0]\, \, \exists \, \text{ an } \alpha-\text{chain }
u_1 < \dots < u_m \notag \\
&\text{ supported by } L \text{ with }
\|\sum_{i=1}^m u_i \| \geq c \, m \} \notag 
\end{align}
we have that \(C \leq \tau \leq 1\).
\end{Lem}
\begin{proof}
We recall the following fact about Schreier families established in \cite{otw}:
Given \(\xi < \eta < \omega_1\), there exists some \(n \in \mathbb{N}\) 
such that for every \(F \in S_\xi\) with \(\min F \geq n\) we have that
\(F \in S_\eta\). Repeated applications of this fact now yield
\(N_0 \in [N]\) with the following property: Let \(j \in N_0\) 
and let \(F \subset N_0\) with \(j < \min F\). Suppose that \(F \in S_{\alpha_i}\) for some
\(i \in \mathbb{N}\) with
\(i \leq j\). Then, \(F \in S_{\alpha_j}\).

Next, let \(A\) denote the set of all \(c\)'s which appears in the definition of
\(\tau \). We need to show that \(C \in A\). Our assumptions yield that for every 
\(L \in [N]\), every \(\epsilon > 0\) and all \(\beta < \alpha\), there exists a 
\((\beta, \epsilon)\) block supported by \(L\). It follows now from this and 
the selection of \(N_0\), that every \(L \in [N_0]\) supports an \(\alpha\)-chain. 
We infer now from the above and the manner \(C\) is defined, that \(C \in A\) and so
the lemma is proved.
\end{proof}
Our next two lemmas will help us to produce an \((\alpha, \epsilon)\) block
from the members of an \(\alpha\)-chain. 
\begin{Lem} \label{L3}
Under the assumptions of Lemma \ref{L2}, 
let \(N_0 \in [N]\) and \(\tau > 0\) be as in the conclusion of that lemma.
Then, for every \(\delta > 0\) there exists \(M_0 \in [N_0]\) 
fulfilling the following property: For every \(\alpha\)-chain
\(u_1 < \dots < u_n\) supported by \(M_0\) and for every choice
\(t_2, \dots, t_n\) of elements of \(K\) with \(t_i\) strongly norming \(u_i\) for
\(i= 2, \dots, n\), we have that
\[ \|u | J \| \leq \delta + (\tau + 2 \delta) 
\sum_{i=2}^n (u|J \cap \mathrm{ supp } \, u_i )(t_i)\]
for all \(J \subset \cup_{i=1}^n \mathrm{ supp } \, u_i \), where, in the above,
we have set \(u = \frac{\sum_{i=1}^n u_i}{\| \sum_{i=1}^n u_i \|}\).
\end{Lem}
\begin{proof}
The definition of \(\tau\) allows us to choose \(M_1 \in [N_0]\) such that
\(\|\sum_{i=1}^n u_i \|\) \(<\) \((\tau + \delta) n\), for every \(\alpha\)-chain
\(u_1 < \dots < u_n\) supported by \(M_1\). Next, let \(m \in M_1\) satisfy
\begin{equation} \label{EL3}
(\tau + 2 \delta)(m-1 - \sum_{s=2}^m \epsilon_s) > (\tau + \delta)m.
\end{equation}
Now let \(u_1 < \dots < u_n\) be an \(\alpha\)-chain supported by \(M_1\) with \(m < n\),
and assume that \(t_i \in K\) strongly norms \(u_i\) for \(2 \leq i \leq n\).
Given \(t \in K\) define
\[I_i^t = \{j \in \mathrm{ supp } \, u_i : \, f_j(t) \geq
(\tau + 2 \delta) f_j(t_i) \}, \text{ for } i= 2, \dots, n.\]
We also set \(I^t = \{i \in \{2, \dots, n\}: \, \|u_i|I_i^t \| \geq \epsilon_{d_{i-1}}^{1/2}\}\)
(recall that \(d_i = \max \mathrm{ supp } \, u_i\)).

{\bf Claim}: \(|I^t| \leq 2m\), for all \(t \in K\).

Once this claim is established, the proof of the lemma is completed as follows:
Let \(m \in M_1\) satisfy \eqref{EL3} and choose \(m_0 \in M_1\) such that
\begin{equation} \label{EL32}
2m + 1 + \sum_{s \geq m_0} \epsilon_s^{1/2} < \delta C m_0.
\end{equation}
Put \(M_0 = \{j \in M_1 : \, j > m_0\}\).
Let \(u_1 < \dots < u_n\) be an \(\alpha\)-chain supported by \(M_0\).
Put \(J_i = \mathrm{ supp } \, u_i\) and write \(u_i = \sum_{j \in J_i} \lambda_j f_j\),
where each \(\lambda_j\) is positive, for all \(i \leq n\).
Note also that
\begin{equation} \label{EL33}
D = \|\sum_{i=1}^n u_i \| \geq Cn > C m_0.
\end{equation}
Assume that \(t_i \in K\) strongly norms \(u_i\) for \(2 \leq i \leq n\).
Given \(t \in K\), let \(I_i^t\) (\(2 \leq i \leq n\)) and \(I^t\) be defined as in the preceding
paragraph. Let \(J \subset \cup_{i=1}^n J_i\).
Suppose first that \(i \in I^t\). We have the following estimates:
\begin{align}
(u_i |J)(t) &= (u_i | J \cap I_i^t)(t) + (u_i | J \setminus I_i^t)(t) \notag \\
 &\leq 1 + (\tau + 2 \delta) \sum_{j \in J_i \cap J} \lambda_j f_j(t_i). \notag
\end{align}
Our claim now implies that
\begin{equation} \label{EL34}
\sum_{i \in I^t} (u_i | J)(t) \leq 2m + 
(\tau + 2 \delta) \sum_{i \in I^t} \sum_{j \in J_i \cap J} \lambda_j f_j(t_i). 
\end{equation}
When \(i \geq 2 \) but \(i \notin I^t\), we have that
\begin{align}
(u_i |J)(t) &= (u_i | J \cap I_i^t)(t) + (u_i | J \setminus I_i^t)(t) \notag \\
&\leq \epsilon_{d_{i-1}}^{1/2} + (\tau + 2 \delta) \sum_{j \in J_i \cap J} \lambda_j f_j(t_i). \notag
\end{align}
We thus obtain that
\begin{equation} \label{EL35}
\sum_{i \geq 2, \, i \notin I^t} (u_i |J)(t) \leq \sum_{i=2}^n \epsilon_{d_{i-1}}^{1/2}
+ (\tau + 2 \delta) \sum_{i \geq 2, \, i \notin I^t} 
\sum_{j \in J_i \cap J} \lambda_j f_j(t_i). 
\end{equation}
Equations \eqref{EL34} and \eqref{EL35} lead us to
\begin{equation} \label{EL36}
\sum_{i=2}^n (u_i | J)(t) \leq 2m + \sum_{i=2}^n \epsilon_{d_{i-1}}^{1/2} +
(\tau + 2 \delta) \sum_{i=2}^n 
\sum_{j \in J_i \cap J} \lambda_j f_j(t_i). 
\end{equation}
It follows now from \eqref{EL36} that
\begin{align}
(u| J) (t) &\leq (1/D) (2m +1 + \sum_{i=2}^n \epsilon_{d_{i-1}}^{1/2})
+ (1/D) (\tau + 2 \delta) \sum_{i=2}^n 
\sum_{j \in J_i \cap J} \lambda_j f_j(t_i) \notag \\
&< (1/(Cm_0))(2m +1 + \sum_{s \geq m_0} \epsilon_s^{1/2}) +
(\tau + 2 \delta) \sum_{i=2}^n (u |J \cap J_i)(t_i), \notag \\
&\text{ by } \eqref{EL33}. \notag
\end{align}
Taking in account \eqref{EL32}, we deduce from the above that
\[(u | J)(t) < \delta + (\tau + 2 \delta) \sum_{i=2}^n (u | J \cap J_i)(t_i),\]
for all \(t \in K\) and all \(J \subset \cup_{i=1}^n J_i\), as desired.

We next give the proof of the claim. Note that by replacing \(M_1\) by a suitable infinite subset
if necessary, there will be no loss of generality in assuming that 
\(\epsilon_j < \epsilon_i^2\) for all \(i < j\) in \(M_1\).
Arguing as we did when selecting \(N_0\) in the proof of Lemma \ref{L2},
we may also assume that \(M_1\) enjoys the following property:
Let \(j \in M_1\) 
and let \(F \subset M_1\) with \(j < \min F\). Suppose that \(F \in S_{\alpha_i}\) for some
\(i \in \mathbb{N}\) with
\(i \leq j\). Then, \(F \in S_{\alpha_j}\).

Assume, on the contrary, that the claim is false.
Then for some \(t \in K\) we would have that \(|I^t| > 2m\).
We may thus choose indices \(2 \leq i_1 < \dots < i_{2m} \leq n\)
such that \(\|u_{i_k} | I_{i_k}^t \| \geq \epsilon_{d_{i_k - 1}}^{1/2}\)
for all \(k=1, \dots , 2m\). If we let \(v_k = \frac{u_{i_k} | I_{i_k}^t}{\|u_{i_k} | I_{i_k}^t \|}\)
for all \(k \leq 2m\), then we deduce from Lemma \ref{L1} that \(v_k\)
is an \((\alpha_{d_{i_k -1}}, \epsilon_{d_{i_k -1}}^{1/2})\) block for all
\(k \leq 2m\). Our assumptions on \(M_1\) now yield that
\(f_m < v_4 < v_6 < \dots < v_{2m}\) is an \(\alpha\)-chain supported by \(M_1\).
Set \(l_1 = m\) and \(l_k = \max \mathrm{ supp } \, v_{2k}\) for \(k=2, \dots, m\).
We observe that for all \(k = 2, \dots , m\),
\[1 = \|v_{2k}\| \leq (1 + \epsilon_{d_{i_{2k} -1 }}^{1/2} ) v_{2k}(t_{i_{2k}})
< (1 + \epsilon_{l_{k-1}}) v_{2k}(t_{i_{2k}}), \]
by our assumptions on \(M_1\) and since \(t_{i_{2k}}\) strongly norms \(v_{2k}\) for
all \(k \leq m\), thanks to Lemma \ref{L1}.
We infer from this observation that
\begin{align} 
v_{2k}(t) &= (\|u_{i_{2k}} | I_{i_{2k}} \|)^{-1} \sum_{j \in I_{i_{2k}}^t}
\lambda_j f_j (t) \label{EL37} \\ 
&\geq (\|u_{i_{2k}} | I_{i_{2k}} \|)^{-1} (\tau + 2 \delta) 
\sum_{j \in I_{i_{2k}}^t}
\lambda_j f_j (t_{i_{2k}}) \notag \\
&= (\tau + 2 \delta) v_{2k}(t_{i_{2k}}) \notag \\
&\geq (\tau + 2 \delta) (1 + \epsilon_{l_{k-1}})^{-1}, \forall \, k=2, \dots , m. \notag
\end{align}
We finally have the estimate
\begin{align}
\|f_m + \sum_{k=2}^m v_{2k} \| &\geq \sum_{k=2}^m v_{2k}(t) \notag \\
&\geq (\tau + 2 \delta) \sum_{k=2}^m (1 + \epsilon_{l_{k-1}})^{-1}, \text{ by }
\eqref{EL37} \notag \\
&= (\tau + 2 \delta) \sum_{k=2}^m (1 - \frac{\epsilon_{l_{k-1}}}{1 + \epsilon_{l_{k-1}}})
\geq (\tau + 2 \delta) ( m-1 - \sum_{s=2}^m \epsilon_s) \notag \\
&> (\tau + \delta) m, \text{ by } \eqref{EL3}, \notag
\end{align}
which contradicts the choice of \(M_1\), as \(f_m < v_4 < v_6 < \dots < v_{2m}\)
is an \(\alpha\)-chain supported by \(M_1\).
This proves the claim and completes the proof of the lemma.
\end{proof}
\begin{Lem} \label{L4}
Under the assumptions of Lemma \ref{L2}, 
let \(N_0 \in [N]\) and \(\tau > 0\) be as in the conclusion of that lemma.
Then for every \(M \in [N_0]\) and \(0 < \delta < \tau /2\),
there exists an \(\alpha\)-chain
\(u_1 < \dots < u_n\)  supported by \(M\) and satisfying the following property:
There exist \(t_0 \in K\) and \(t_2, \dots, t_n\) in \(K\)
with \(t_i\) strongly norming \(u_i\) for all \(i=2, \dots, n\)
and such that 
\[f_j(t_0) \geq (\tau - 2 \delta) f_j(t_i),\]
for all \(j \in \mathrm{ supp } \, u_i\) and
\(i=2, \dots, n\).
\end{Lem}
\begin{proof}
Let \(M \in [N_0]\). We may assume without loss of generality that
\(\epsilon_j < \epsilon_i^2\) for all \(i < j\) in \(M\).
Arguing as we did when selecting \(N_0\) in the proof of Lemma \ref{L2},
we may assume in addition to the above, that \(M\) posesses the following property:
If \(j \in M\) 
and \(F \subset M\) with \(j < \min F\) is so that \(F \in S_{\alpha_i}\) for some
\(i \in \mathbb{N}\) with
\(i \leq j\), then \(F \in S_{\alpha_j}\).

We first choose some \(n \in M\) with \(n > 2\) and then choose \(M_1 \in [M]\)
with \( 2n + 2 < (\delta /2) \min M_1\). By the definition of \(\tau\), we can select
an \(\alpha\)-chain \(w_1 < \dots < w_p\) supported by \(M_1\) and such that
\(\|\sum_{i=1}^p w_i \| \geq (\tau - \delta)p\). Choose \(t_0 \in K\) such that
\begin{equation} \label{EL41}
\sum_{i=1}^p w_i(t_0) = \|\sum_{i=1}^p w_i \| \geq (\tau - \delta) p.
\end{equation}
Put \(J_i = \mathrm{ supp } \, w_i\) and write \(w_i = \sum_{j \in J_i} \lambda_j f_j\),
where each \(\lambda_j\) is positive, for all \(i \leq p\).
We define 
\[I_i = \{ j \in J_i: \, f_j(t_0) \geq (\tau - 2 \delta) f_j(x_i) \},\]
where \(x_i \in K \) strongly norms \(w_i\) for \(i = 2, \dots, p\).
Set
\(I = \{i \in \{2, \dots, p\}: \, \|w_i|I_i \| \geq \epsilon_{d_{i-1}}^{1/2}\}\)
(recall that \(d_i = \max \mathrm{ supp } \, w_i\) for all \(i \leq p\)).
We claim that \(I\) contains more than \(2n\) elements.
Indeed, if not, then 
\begin{align}
\sum_{i=1}^p w_i(t_0) &= 1 + \sum_{i \in I} w_i (t_0) + 
\sum_{i \geq 2, \, i \notin I} w_i (t_0) \notag \\
&= 1 + \sum_{i \in I} w_i(t_0) + \sum_{i \geq 2, \, i \notin I}
[(w_i | I_i)(t_0) + (w_i | J_i \setminus I_i)(t_0)] \notag \\
&\leq 1 + 2n + \sum_{i=2}^p
\epsilon_{d_{i-1}}^{1/2} + 
\sum_{i \geq 2, \, i \notin I} \sum_{j \in J_i \setminus I_i}
\lambda_j f_j(t_0), 
\text{ since } |I| \leq 2n \notag \\
&\leq 1 + 2n + \sum_{i=2}^p
\epsilon_{d_{i-1}}^{1/2} +
(\tau - 2 \delta) \sum_{i \geq 2, \, i \notin I} \sum_{j \in J_i \setminus I_i}
\lambda_j f_j(x_i) \notag \\
&\leq  2 + 2n + (\tau - 2 \delta)p. \notag
\end{align}
Combining this with \eqref{EL41} we obtain that
\[(\tau - \delta)p \leq 2 + 2n + (\tau - 2 \delta)p\]
and so 
\[\delta p \leq 2 + 2n < (\delta /2) \min M_1 \leq (\delta /2 )p, \]
as \(p \in M_1\). This contradiction proves the claim.
We may now choose indices \(2 \leq i_1 < \dots < i_{2n} \leq p\)
such that \(\|w_{i_k} | I_{i_k} \| \geq \epsilon_{d_{i_k - 1}}^{1/2}\)
for all \(k=1, \dots , 2n\). Set \(v_k = (\|w_{i_k} | I_{i_k} \|)^{-1} (w_{i_k} | I_{i_k})\),
for all \(k \leq 2n\).   
We infer now from Lemma \ref{L1} that \(v_k\)
is an \((\alpha_{d_{i_k -1}}, \epsilon_{d_{i_k -1}}^{1/2})\) block strongly normed
by \(x_{i_k}\) for all
\(k \leq 2n\). Our assumptions on \(M\) now yield that
\(f_n < v_4 < v_6 < \dots < v_{2n}\) is an \(\alpha\)-chain supported by \(M\)
for which \(f_j(t_0) \geq (\tau - 2 \delta) f_j(x_{i_{2k}})\) for all
\(j \in \mathrm{ supp } \, v_{2k}\) and all \(k =2, \dots , n \).
Put \(u_1 = f_n\) and \(u_k = v_{2k}\) for all \(k \leq n\).
Then \(u_1 < \dots < u_n\) is the desired \(\alpha\)-chain, since \(t_k = x_{i_{2k}}\) strongly norms
\(u_k\) for all \(k =2, \dots, n\).
\end{proof}
\begin{proof}[Proof of Theorem \ref{T3}]
We use transfinite induction on \(\alpha < \omega_1\).
The case \(\alpha = 0\) was settled in Lemma \ref{L0}.
Assume that \(\alpha \geq 1\) and that the assertion of the theorem holds
for all ordinals \(\beta < \alpha\). Choose \(N_0 \in [N]\) and \(\tau > 0\)
according to Lemma \ref{L2}. Let \(0 < \epsilon < 1 \)
and choose
\(0 < \delta < \tau / 2\) and \(0 < \epsilon_0 < \epsilon\) satisfying
\[ \delta < \epsilon_0^2 \, \text{ and } \, 
\frac{\tau + 2 \delta}{(1 - \epsilon_0)(\tau - 2 \delta)} < 1 + \epsilon.\]
Let \(M_0 \in [N_0]\) satisfy the conclusion of Lemma \ref{L3}.
We apply Lemma \ref{L4} to find an \(\alpha\)-chain
\(u_1 < \dots < u_n\) supported by \(M_0\) for which there exist
\(t_0 \) and \(t_2, \dots , t_n \) in \(K\) with \(t_i\) strongly norming
\(u_i\) for all \(i=2, \dots, n\), and such that
\begin{equation} \label{ET21}
f_j(t_0) \geq (\tau - 2 \delta) f_j (t_i), \, \forall \, j \in \mathrm{ supp } \,
u_i, \, \forall \, i=2, \dots, n.
\end{equation}
We set \(u = \frac{\sum_{i=1}^n u_i}{\|\sum_{i=1}^n u_i \|}\).
We are going to show that \(u\) is an \((\alpha, \epsilon)\) block.

We first consider some \(J \subset \mathrm{ supp } \, u\) such that
\(\|u | J\| \geq \epsilon\).
Since \(u_1 < \dots < u_n\) is an \(\alpha\)-chain supported 
by \(M_0\), Lemma \ref{L3} yields
\begin{align}
\|u | J \| &\leq \delta + (\tau + 2 \delta) \sum_{i=2}^n (u | J \cap \mathrm{ supp } \, u_i)(t_i) \notag \\
&\leq \delta + \frac{\tau + 2 \delta}{\tau - 2 \delta} 
\sum_{i=2}^n (u | J \cap \mathrm{ supp } \, u_i)(t_0), \text{ by } \eqref{ET21}. \notag
\end{align}
It follows now that
\begin{align}
\|u | J \| &\leq \delta + \frac{\tau + 2 \delta}{\tau - 2 \delta} (u | J)(t_0) 
< \epsilon_0^2 + \frac{\tau + 2 \delta}{\tau - 2 \delta} (u | J)(t_0) \notag \\
&< \epsilon_0 \|u | J\| + \frac{\tau + 2 \delta}{\tau - 2 \delta} (u | J ) (t_0),
\text{ as } \|u | J \| \geq \epsilon, \notag
\end{align}
whence \(\|u | J \| < 
\frac{\tau + 2 \delta}{(1 - \epsilon_0)(\tau - 2 \delta)} (u | J)(t_0) \).
We deduce now from our initial choices, that
\[ \|u | J \| < ( 1 + \epsilon) ( u | J )(t_0)\]
and thus \(t_0\) strongly norms \(u\).

We finally show that \(\|u | J \| < \epsilon^2\) for every \(J \subset \mathrm{ supp } \, u\)
with \(J \in S_\alpha\). Set \(D = \| \sum_{i=1}^n u_i \|\).
We can certainly assume, without loss of generality, that
\[\sum_{l \in N} l \epsilon_l^2 < \epsilon^2/2, \text{ and that } \,
\min N > 2 / (C \epsilon^2).\]
We consider now some \(J \subset \mathrm{ supp } \, u\)
with \(J \in S_\alpha\). 
We let \(i_0\) denote the smallest \(i \leq n\) for which
\(J \cap \mathrm{ supp } \, u_i \ne \emptyset\). It follows
that \(\min J \leq d_{i_0}\) (recall that \(d_i = \max \mathrm{ supp } \, u_i\))
and thus, by the definition of Schreier families, we have that
\(J \in S_{\alpha_j + 1}\) for some \(j \leq \min J\). We can therefore
write \(J = \cup_{s=1}^p J_s\) where \(J_1 < \dots < J_p\) are members of \(S_{\alpha_j}\)
and \(p \leq \min J\). We obtain now, from the fact that \(u_1 < \dots < u_n\)
is an \(\alpha\)-chain and because \(j \leq d_{i_0}\), that
\[\|u_i | J_s \| < \epsilon_{d_{i-1}}^2, \, \forall \, s \leq p, \, \forall \,
i = i_0 +1 , \dots, n.\]
Hence, \(\|u_i | J \| \leq p \epsilon_{d_{i-1}}^2 \leq d_{i_0} \epsilon_{d_{i-1}}^2\),
for all \(i = i_0 +1 , \dots , n\).
Summarizing all the above, we get to
\[ \sum_{i = i_0 +1}^n \| u_i | J \| \leq \sum_{i=i_0}^{n-1} d_i \epsilon_{d_i}^2 < \epsilon^2 /2.\]
Because we assumed that \(D \geq Cn > 2/\epsilon^2\),
we also get that 
\((1/D) \|u_{i_0} | J \|  < \epsilon^2 /2 \) and so
\(\|u | J \| < \epsilon^2\), as desired. This completes the proof
of the theorem.
\end{proof}
\section{Proof of the main result} \label{S4}
In this section we present the proof of the main result of this paper.
In what follows, \(K\) is a compact metrizable space, \(X\) is a subspace of \(C(K)\)
which is asymptotic \(\ell_1\) with respect to its normalized Schauder basis \((e_n)\).
We let \(f_n = |e_n|\), for all \(n \in \mathbb{N}\). We also fix a decreasing sequence of
positive scalars \((\epsilon_n)\) with \(\sum_n \epsilon_n^{1/2} < 1\).
Our next proposition uses an idea from Lemma 7 of \cite{b}.
\begin{Prop} \label{mP}
Let \(\mathcal{M} \subset B_{C(K)^*}\) be a \(w^*\)-compact set consisting of positive measures
on \(K\). Suppose that \(\mathcal{M}\) \(\rho\)-norms \(X\) for some \(\rho > 0\). 
Then for every \(\alpha < \omega_1\)
and all \(N \in [\mathbb{N}]\), \(N = (n_i)\), there exist \(I \in [\mathbb{N}]^{< \infty}\)
and \(\mu \in \mathcal{M}\) satisfying the following:
\begin{enumerate}
\item \(\mu ([ |f_{n_{2i}}| \geq \epsilon_{n_{2i-1}}]) \geq \rho^2 / 5\), for all \(i \in I\).
\item \(\{n_{2i}: \, i \in I\} \notin S_\alpha\).
\end{enumerate}
\end{Prop}
\begin{proof}
We first choose \(0 < \epsilon < \rho / 2\) such that
\[\biggl (\frac{\rho - 2 \epsilon}{2 + \epsilon} \biggr )^2 > \epsilon^2 + \rho^2 / 5 . \]
We set \(D = \frac{2 + \epsilon}{\rho - 2 \epsilon}\).
We may assume, without loss of generality, that \(\sum_i \epsilon_{n_{2i-1}}\)
\( <\) \(\epsilon\).

Theorem \ref{T3} enables us to find an \((\alpha, \epsilon)\) block
\(u = \sum_i a_i f_{n_{2i}}\) supported by \(\{n_{2i}: \, i \in \mathbb{N}\}\).
Let \(t_0 \in K\) strongly norm \(u\). We thus have that
\begin{align} 
&\|\sum_{i \in J} a_i f_{n_{2i}} \| \leq (1 + \epsilon) 
 \sum_{i \in J} a_i f_{n_{2i}}(t_0), \, \forall \, 
J \subset \mathbb{N} \text{ with }
\|\sum_{i \in J} a_i f_{n_{2i}} \| \geq \epsilon \label{EmP1} \\ 
&\text{ and } \|\sum_{i \in J} a_i f_{n_{2i}} \| < \epsilon^2, \, \forall \,
J \subset \mathbb{N} \text{ with } \{n_{2i}: \, i \in J \} \in S_\alpha. \notag
\end{align}
We can select \(\mu \in \mathcal{M}\) such that
\(\int_K \, (\sum_i a_i f_{n_{2i}} ) \, d \mu \geq \rho\).
To see this, first choose \(t_1 \in K\) with
\(\sum_i a_i f_{n_{2i}}(t_1) = 1\). We can now find signs \(\sigma_i\) so that
\( \sum_i a_i \sigma_i e_{n_{2i}}(t_1) = 1\). It follows that
\(\| \sum_i a_i \sigma_i e_{n_{2i}}\| = 1\). Since \(\mathcal{M}\) \(\rho\)-norms \(X\),
there exists \(\mu \in \mathcal{M}\) such that
\(|\int_K \, (\sum_i a_i \sigma_i e_{n_{2i}}) \, d \mu | \geq \rho\)
and hence \(\mu\) does the job.

We next put 
\[\delta_i = \int_K \, f_{n_{2i}} \, d \mu , \forall \, i \in \mathbb{N}, \text{ and, }
I_0 = \{i \in \mathbb{N}: \, n_{2i} \in \mathrm{ supp } \, u, \text{ and, } \delta_i > 0 \}.\]
We also set \(\phi_i = \chi_{[f_{n_{2i}} \geq D \delta_i]}\), for all \(i \in I_0\).
Chebyshev's inequality yields 
\begin{equation} \label{EmP2}
\int_K \, \phi_i \, d \mu \leq 1 / D, \, \forall \, i \in I_0.
\end{equation}
We now have the following estimates:
\begin{align}
\int_K \, (\sum_i a_i f_{n_{2i}} ) \, d \mu &= \int_K \, (\sum_{i \in I_0} a_i f_{n_{2i}} ) \, d \mu
\label{EmP3} \\
&= \int_K \, (\sum_{i \in I_0} a_i f_{n_{2i}} \phi_i) \, d \mu +
\int_K \, (\sum_{i \in I_0} a_i f_{n_{2i}} (1 - \phi_i)) \, d \mu \notag \\
&= \int_K \, (\sum_{i \in I_0} a_i f_{n_{2i}} \phi_i ) \, d \mu +
\sum_{i \in I_0} a_i \int_{ [ f_{n_{2i}} < D \delta_i]} \, f_{n_{2i}} \, d \mu \notag \\
&\leq \int_K \, \|\sum_{i \in I_0} a_i \phi_i(t) f_{n_{2i}}  \| \, d \mu (t)
+ \sum_{i \in I_0} a_i \int_{ [ f_{n_{2i}} < D \delta_i]} \, f_{n_{2i}} \, d \mu . \notag
\end{align} 
Let \(K_1 = \{ t \in K : \, \|\sum_{i \in I_0} a_i \phi_i(t) f_{n_{2i}}  \| \geq \epsilon \}\).
Clearly, \(K_1\) is a closed subset of \(K\). Since \(u\) is an \((\alpha, \epsilon)\) block, 
\eqref{EmP1} implies that
\[\|\sum_{i \in I_0} a_i \phi_i(t) f_{n_{2i}}  \| \leq (1 + \epsilon) 
\sum_{i \in I_0} a_i \phi_i(t) f_{n_{2i}}(t_0), \, \forall \, t \in K_1\]
and so
\begin{align} 
\int_{K_1} \, \|\sum_{i \in I_0} a_i \phi_i(t) f_{n_{2i}}  \| \, d \mu(t) &\leq
(1 + \epsilon) \int_{K_1} \, \sum_{i \in I_0} a_i \phi_i(t) f_{n_{2i}}(t_0)  
\label{EmP4} \\
&= (1 + \epsilon) \sum_{i \in I_0} a_i f_{n_{2i}}(t_0) \int_{K_1} \, \phi_i \, d \mu \notag \\
&\leq \frac{1 + \epsilon}{D} \sum_{i \in I_0} a_i f_{n_{2i}}(t_0) 
 \leq (1 + \epsilon)/ D, \text{ by }
\eqref{EmP2}. \notag 
\end{align}  
Note also that \(\int_{K \setminus K_1} \, 
\|\sum_{i \in I_0} a_i \phi_i(t) f_{n_{2i}}  \| \, d \mu(t) \leq \epsilon \) 
and hence we deduce from \eqref{EmP4} that
\begin{equation} \label{EmP5}
\int_K \, \|\sum_{i \in I_0} a_i \phi_i(t) f_{n_{2i}}  \| \, d \mu(t) \leq 
\epsilon + (1 + \epsilon)/D.
\end{equation}
Taking in account \eqref{EmP5}, \eqref{EmP3} gives us that
\begin{align}
\rho &\leq \int_K \, (\sum_i a_i f_{n_{2i}} ) \, d \mu \label{EmP6} \\
&\leq \epsilon + (1 + \epsilon)/D +
\sum_{i \in I_0} a_i \int_{ [ f_{n_{2i}} < D \delta_i]} \, f_{n_{2i}} \, d \mu \notag \\
&\leq \epsilon + (1 + \epsilon)/ D + 
\sum_{i \in I_0} a_i \bigl [\int_{ [ \epsilon_{n_{2i-1}} \leq f_{n_{2i}} < D \delta_i]} \, f_{n_{2i}} \, d \mu 
+ \int_{ [ f_{n_{2i}} < \epsilon_{n_{2i-1}}]} \, f_{n_{2i}} \, d \mu \bigr ] \notag \\
&\leq \epsilon + (1 + \epsilon)/ D +
\sum_{i \in I_0} a_i D \delta_i \mu ( [f_{n_{2i}} \geq \epsilon_{n_{2i-1}}])
+ \sum_{i \in I_0} a_i \epsilon_{n_{2i-1}} \notag \\
&\leq 2 \epsilon + (1 + \epsilon)/D + D \sum_i a_i \delta_i  
\mu ( [f_{n_{2i}} \geq \epsilon_{n_{2i-1}}]), \text{ as }
\sum_i \epsilon_{n_{2i - 1}} < \epsilon. \notag
\end{align}
A straightforward computation shows that
\[ \rho - 2 \epsilon - (1 + \epsilon) / D = D 
\biggl (\frac{\rho - 2 \epsilon}{2 + \epsilon} \biggr )^2\]
and so \eqref{EmP6} leads to 
\begin{equation} \label{EmP7}
\sum_i a_i \delta_i  
\mu ( [f_{n_{2i}} \geq \epsilon_{n_{2i-1}}]) \geq 
\biggl (\frac{\rho - 2 \epsilon}{2 + \epsilon} \biggr )^2.
\end{equation}
We now define \(I = \{ i \in \mathbb{N}: \, n_{2i} \in \mathrm{ supp } \, u, \text{ and, }
\mu([f_{n_{2i}} \geq \epsilon_{n_{2i-1}}]) \geq \rho^2 /5\}\).
We claim that \(\{n_{2i} : \, i \in I \} \notin S_\alpha\).
Indeed, suppose the claim is false. Then, since \(u\) is an \((\alpha, \epsilon)\) block,
\eqref{EmP1} yields
\(\|\sum_{i \in I} a_i f_{n_{2i}} \| < \epsilon^2\) and thus
\[\sum_{i \in I} a_i \delta_i \mu( [f_{n_{2i}} \geq \epsilon_{n_{2i-1}}]) 
\leq \sum_{i \in I} a_i \delta_i = \int_K \, \sum_{i \in I} a_i f_{n_{2i}} \, d \mu
< \epsilon^2.\] 
We also observe that \(\sum_i a_i \delta_i = \int_K \, u \, d \mu \leq 1\)
and therefore,
\[\sum_{i \notin I} a_i \delta_i \mu( [f_{n_{2i}} \geq \epsilon_{n_{2i-1}}]) 
\leq \rho^2 /5.\]
We finally infer from the above and \eqref{EmP7} that
\[\biggl (\frac{\rho - 2 \epsilon}{2 + \epsilon} \biggr )^2 < \epsilon^2 + \rho^2/ 5,\]
which contradicts our initial choice of \(\epsilon\). Thus, our claim holds and 
so \(\mu \in \mathcal{M}\) and \(I \in [\mathbb{N}]^{< \infty}\) are desired.
The proof of the proposition is now complete.
\end{proof}
\begin{Cor} \label{C1}
Let \(\mathcal{M} \subset B_{C(K)^*}\) be a \(w^*\)-compact set consisting of positive measures
on \(K\). Suppose that \(\mathcal{M}\) \(\rho\)-norms \(X\) for some \(\rho > 0\). 
Then for every \(M \in [\mathbb{N}]\), \(M = (m_i)\), the family
\begin{align}
\mathcal{F}_{M} = \bigl \{ &\{m_{2i}: \, i \in I \} : \,
I \in [\mathbb{N}]^{< \infty}, \text{ and, } \exists \, \mu \in \mathcal{M} \text{ with } \notag \\
&\mu ( [f_{m_{2i}} \geq \epsilon_{m_{2i-1}}]) \geq \rho^2/5, \, \forall \, i \in I \bigr \} \notag
\end{align}
is not compact in the topology of pointwise convergence.
\end{Cor}
\begin{proof}
Suppose the assertion of the corollary is false. Then \(\mathcal{F}_M\) would be
pointwise compact for some \(M \in [\mathbb{N}]\), \(M= (m_i)\).
It follows by the Mazurkiewicz-Sierpinski theorem \cite{ms}, that
\(\mathcal{F}_M\) is homeomorphic to the ordinal interval
\([1, \omega^\beta d]\), for some \(\beta < \omega_1\) and \(d \in \mathbb{N}\).
Set \(\alpha = \beta + 1\). It is easily seen that
\(\mathcal{F}_M\) is hereditary and so we can apply the result of \cite{g} (see also
\cite{ju})
to obtain an infinite subset
\(L\) of \(\{m_{2i}: \, i \in \mathbb{N}\}\) so that
\(\mathcal{F}_{M}[L] \subset S_\alpha\). This is possible since by \cite{aa}
\(S_\alpha[P]\) is homeomorphic to \([1, \omega^{\omega^\alpha}]\), for all \(P \in [\mathbb{N}]\),
which in turn, is not homeomorphic to a subset of \([1, \omega^\beta d]\).
Set \(N = L \cup \{m_{2i-1}: \, i \in \mathbb{N}, \text{ and, } m_{2i} \in L \}\).
We infer from Proposition \ref{mP} that there exist \(I \in [\mathbb{N}]^{< \infty}\)
and \(\mu \in \mathcal{M}\) such that
\(\mu([f_{n_{2i}} \geq \epsilon_{n_{2i-1}}]) \geq \rho^2 /5 \) for all \(i \in I\) and
\(\{n_{2i} : \, i \in I \} \notin S_\alpha \). However, 
\(\{n_{2i} : \, i \in I \} \in \mathcal{F}_M[L]\) \(\subset S_\alpha\). This contradiction
proves the corollary.
\end{proof}
\begin{Cor} \label{C2}
For every \(M \in [\mathbb{N}]\) there exists \(L \in [M]\), \(L = (l_i)\), with the property
that for every \(n \in \mathbb{N}\) there exists a \(\mu \in \mathcal{M}\) so that
\[\mu( [ f_{l_{2i}} \geq \epsilon_{l_{2i-1}}]) \geq \rho^2 /5\] for every \(i \leq n\).
\end{Cor}
\begin{proof}
We know, thanks to Corollary \ref{C1}, that \(\mathcal{F}_M\) is hereditary but not pointwise compact.
We deduce from this, that there exists \(J \in [\mathbb{N}]\), \(J = (j_i)\), such that
\(\{m_{2j_1}, \dots, m_{2j_n}\} \in \mathcal{F}_M \) for all \(n \in \mathbb{N}\).
We need only take \(l_{2i} = m_{2j_i}\) and \(l_{2i-1} = m_{2j_i -1 }\) for all
\(i \in \mathbb{N}\) to produce the required \(L \in [M]\).
\end{proof}
\begin{Cor} \label{C3}
Suppose that \(K\) is totally disconnected and let \(\mathcal{M}\) be a 
\(w^*\)-compact subset of \(B_{C(K)^{*}}\) which norms \(X\).
Then \(\mathcal{M}\) is not separable in norm.
\end{Cor}
\begin{proof}
Assume on the contrary, that \(\mathcal{M}\) is separable in norm.
Note first that a result due to R. Haydon \cite{h} implies that \(X^*\) is separable.
It follows from this and Rosenthal's theorem \cite{r2},
that there is no loss of generality in assuming that \((e_n)\)
is a normalized, shrinking Schauder basis for \(X\). In particular, \((f_n)\) is a 
normalized, weakly null sequence in \(C(K)\).

We next observe that \(X\), being asymptotic \(\ell_1\), can not contain an isomorph of \(c_0\)
and subsequently, according to a result of Bourgain (Lemma 4 in \cite{b}), there exists
a \(w^*\)-compact set \(\mathcal{N} \subset B_{C(K)^*}\) of positive measures on \(K\)
which is separable in norm, and such that \(\mathcal{N}\) \(\rho\)-norms \(X\) for some
\(\rho > 0\). Define
\begin{align}
\mathcal{D}= \{L \in [\mathbb{N}], L= (l_i) : \, &\forall \, n \in \mathbb{N}, \,
\exists \, \mu \in \mathcal{N} \text{ with } \notag \\
&\mu ( [ f_{l_{2i}} \geq \epsilon_{l_{2i-1}}]) 
\geq \rho^2 /5, \, \forall \, i \leq n \}.\notag
\end{align}
Clearly, \(\mathcal{D}\) is closed in the topology of pointwise convergence. 
The infinite Ramsey theorem (\cite{ell}, \cite{o1}) and Corollary \ref{C2}
now yield \(N \in [\mathbb{N}]\) satisfying \([N] \subset \mathcal{D}\).
We finally deduce from Theorem \ref{T2} that \(\mathcal{N}\) is not separable
in norm. This contradiction shows that \(\mathcal{M}\) is not separable in norm,
as required. 
\end{proof}

\end{document}